\documentclass[twoside]{amsart}

\usepackage{amssymb}
\usepackage{mathrsfs} 
\usepackage{ifthen}

\newtheorem*{theorem}{Theorem}

\newtheorem{q}{Question}
\theoremstyle{definition}
\newtheorem{example}{Example}

\newcounter{oldbookproblem}
\newcounter{newbookproblem}
\newcommand{\bookproblem}[1][1]{%
  \mbox{}%
  \addtocounter{oldbookproblem}{1}%
  \addtocounter{newbookproblem}{#1}%
  \ifthenelse
    {\equal{#1}{1}}
    {
      \ifthenelse
        {\isodd{\thepage}}
        {\marginpar{\bf\arabic{newbookproblem}\!~?}}
        {\marginpar{\flushright\bf?\!~\arabic{newbookproblem}}}%
    }
    {
       \ifthenelse
        {\isodd{\thepage}}
        {\marginpar{\bf\arabic{oldbookproblem}--\arabic{newbookproblem}\!~?}}
        {\marginpar{\flushright\bf?\!~\arabic{oldbookproblem}--\arabic{newbookproblem}}}
       \addtocounter{oldbookproblem}{-1}%
       \addtocounter{oldbookproblem}{#1}%
    }%
}

\newcommand{\abs}[1]{\lvert#1\rvert}
\def\norm#1{\left\Vert#1\right\Vert}

\def\R {{\mathbb R}}

\def\I {{\mathbb I}}
\def\Q {{\mathbb Q}}
\def\N{{\mathbb N}}
\def\Ur{{\mathbb U}}
\def\e{{\varepsilon}}

\def\Z {{\mathbb Z}}
\def\T {{\mathbb T}}
\def\s{{\mathbb S}}

\def\M{{\mathscr{M}}}

\def\H{{\mathcal H}}

\def\Iso{{\mathrm{Iso}\,}}

\def\RUCB{{\mathrm{RUCB}\,}}

\def\Aut{{\mathrm{Aut}\,}}

\def\Homeo{{\mathrm{Homeo}\,}}

\begin{document}
\title{Forty-plus annotated questions about large topological groups}
\author{Vladimir Pestov}
\begin{abstract}
This is a selection of open problems dealing with ``large'' 
(non locally compact) topological groups and concerning
extreme amenability (fixed point on compacta property), oscillation stability,
universal minimal flows and other aspects of universality, and 
unitary representations.
\end{abstract}
\thanks{Research program by the author has been supported by the
NSERC operating grant (2003--07) and the University of Ottawa internal
grants (2002-04 and 2004-08).}
\maketitle

A topological group $G$ is \textit{extremely amenable,} or has the
\textit{fixed point on compacta property,} if every continuous action of $G$
on a compact Hausdorff space has a $G$-fixed point.
Here are some important examples of such groups.

\begin{example}
\label{ex:one}
The unitary group $U(\ell^2)$ of the separable Hilbert space $\ell^2$
with the strong operator topology (that is, the topology of pointwise convergence
on $\ell^2$) (Gromov and Milman \cite{GrM}).
\end{example}

\begin{example}
\label{ex:gfw}
The group $L^1((0,1),\T)$ of all equivalence classes of Borel
maps from the unit interval to the circle with the $L^1$-metric
$d(f,g) = \int_0^1 \abs{f(x)-g(x)}\,dx$
(Glasner \cite{Gl1}, Furstenberg and Weiss, unpublished).
\end{example}

\begin{example}
\label{ex:homeo}
The group $\Aut(\Q,\leq)$ of all order-preserving
bijections of the rationals, equipped with the natural Polish group
topology of pointwise convergence on $\Q$ considered as a discrete space and,
as an immediate corollary,
the group
$\Homeo_+[0,1]$ of all homeomorphisms of the closed unit interval, preserving
the endpoints, equipped with the compact-open topology
(the present author \cite{P98a}).
\end{example}

The above property is not uncommon among concrete ``large'' 
topological groups coming from diverse parts of mathematics. 
In addition to the above quoted articles, we
recommend \cite{GP2,KPT} and the book \cite{P05}.

The group in example \ref{ex:gfw} 
is monothetic, that is, contains a dense subgroup
isomorphic to the additive group of integers $\Z$. 
Notice that every abelian
extremely amenable group $G$ is {\em minimally almost periodic,}
that is, admits
no non-trivial continuous characters 
(the book \cite{DPS} is a useful reference): 
indeed, if $\chi\colon G\to\T$ is such
a character, then $(g,z)\mapsto \chi(g)z$ defines a continuous action of $G$
on $\T$ without fixed points. The converse remains open.

\begin{q}[Eli Glasner \cite{Gl1}]
\bookproblem%
Does there exist a monothetic topological group that is minimally almost
periodic but not extremely amenable?
\label{q:glasner}
\end{q}
An equivalent question
is: does there exist a topology on the group $\Z$ of integers making it into a topological group that admits a
free action on a compact space but has no non-trivial characters?

Suppose the answer to the above question is in the positive, and let
$\tau$ be a minimally almost periodic Hausdorff group topology on $\Z$
admitting a free continuous action on a compact space $X$. 
Let $x_0\in X$. Find an open neighbourhood $V$ of $x_0$ with $1\cdot V\cap V=
\emptyset$.
It is not difficult to verify that the set $S = \{n\in\Z\colon nx_0\in V\}$ is
{\em relatively dense} in $\Z$, that is, the size of gaps between two
subsequent elements of $S$ is uniformly bounded from above, and at the same time, the closure of $S-S=\{n-m\colon n,m\in S\}$ is a proper subset of $\Z$.
The interior of $S-S$ in the Bohr topology on $\Z$ (the finest
precompact group topology) is therefore not everywhere dense in $(\Z,\tau)$.
Assuming this interior is non-empty,
one can now verify that the $\tau$-closures of elements of the Bohr topology on $\Z$ form a base for a precompact group topology that is nontrivial and coarser than $\tau$, contradicting the assumed minimal almost periodicity of $(\Z,\tau)$. Thus, a positive answer to Glasner's
question would answer in the negative the following very old question from
combinatorial number theory/harmonic analysis, rooted in the classical
work of Bogoliuboff, F\o lner \cite{Fol}, Cotlar and Ricabarra \cite{CotRic}, Veech \cite{V2}, and Ellis and Keynes \cite{EK}:

\begin{q}
\bookproblem%
Let $S$ be a relatively dense subset of the integers. Is the set $S-S$ a
Bohr neighbourhood of zero in $\Z$? 
\end{q}

We refer the reader to Glasner's original work \cite{Gl1} for more on the above. See also \cite{W,P98b,P05}.

\begin{q}
\bookproblem%
Does there exists an abelian minimally almost
periodic topological group acting freely on a compact space?
\end{q}

This does not seem to be equivalent to Glasner's problem, because 
there are examples of minimally almost periodic 
abelian Polish groups whose every monothetic subgroup is discrete, 
such as $L^p(0,1)$ with $0<p<1$.

There are numerous known ways to construct monothetic minimally almost periodic
groups \cite{AHK,DPS,Ba,PZ}. The problem is verifying their (non) extreme amenability.
The most general result presently known asserting non extreme amenability of a
topological group is:

\begin{theorem}[Veech \cite{V}]
Every locally compact group admits a free action on a compact space.
\end{theorem}

Since every locally compact abelian group admits sufficiently many characters,
one cannot employ Veech theorem to answer Glasner's question. Can the result be
extended? Recall that a
topological space $X$ is called a {\em $k_\omega$-space}
(or: a {\em hemicompact} space) if it admits a countable cover $K_n$, 
$n\in\N$ by compact subsets in such a way
that an $A\subseteq X$ is closed if and only if $A\cap K_n$ is closed for
all $n$. For example, every countable $CW$-complex, every second countable
locally compact space, and the free topological group 
\cite{graev} on a compact space are such. 

\begin{q}
\bookproblem%
Is it true that every topological group $G$
that is a $k_\omega$-space admits a free action on a compact space?
\label{q:komega}
\end{q}

\begin{q}
\bookproblem%
Same, for abelian topological groups that are $k_\omega$-spaces.
\end{q}

A positive answer would 
have answered in the affirmative Glasner's question because there are
examples of minimally almost periodic $k_\omega$ group topologies
on the group $\Z$ of integers \cite{PZ}.

Recall that the {\em Urysohn universal metric space} $\Ur$ is
the (unique up to an isometry) 
complete separable metric space that is {\em ultrahomogeneous}
(every isometry between two finite subsets extends to a global self-isometry
of $\Ur$) and {\em universal} ($\Ur$ contains an isometric copy of every
separable metric space) \cite{U25,Ver04,Gr,GK}. The group $\Iso(\Ur)$ of all
self-isometries of $\Ur$, equipped with the topology of pointwise convergence
(which coincides with the compact-open topology), is a Polish topological 
group with a number of remarkable properties. In particular, $\Iso(\Ur)$ is a
universal second-countable topological group \cite{Usp90,Usp98} and is
extremely amenable \cite{P02}.

\begin{q} 
\bookproblem%
Is the group $\Iso(\Ur)$ divisible, that is, does
every element possess roots of every positive natural order?\footnote{Recently Julien Melleray has announced a negative answer (private communication).}
\end{q}

Returning to Glasner's question \ref{q:glasner},
every element $f$ of $\Iso(\Ur)$ generates a monothetic Polish subgroup, 
so one can talk of {\em generic monothetic subgroups} of
$\Iso(\Ur)$ (in the sense of Baire category).

\begin{q}[Glasner and Pestov, 2001, unpublished]
\bookproblem%
Is a generic monothetic subgroup of the isometry group $\Iso(\Ur)$ 
of the Urysohn metric space minimally almost periodic?
\end{q}

\begin{q}[Glasner and Pestov]
\bookproblem%
Is a generic monothetic subgroup of $\Iso(\Ur)$ 
of the Urysohn metric space extremely amenable?
\end{q}

The concept of the universal Urysohn metric space admits numerous
modifications. For instance, one can study 
the universal Urysohn metric space $\Ur_1$ {\em of diameter one}
(it is isometric to 
every sphere of radius $1/2$ in $\Ur$). By analogy with the unitary group
$U(\ell^2)$, it is natural to consider the {\em uniform topology} on the
isometry group $\Iso(\Ur_1)$, given by the bi-invariant uniform metric
$d(f,g)=\sup_{x\in\Ur_1}d_{\Ur_1}(f(x),g(x))$. It is strictly finer than the
strong topology.

\begin{q} 
\bookproblem%
Is the uniform topology on $\Iso(\Ur_1)$ non-discrete?\footnote{According to Julien Melleray (a private communication), the answer is yes.}
\end{q}

\begin{q} 
\bookproblem%
Does $\Iso(\Ur_1)$ possess a uniform neighbourhood of zero covered by
one-parameter subgroups? 
\end{q}

\begin{q} 
\bookproblem%
Does $\Iso(\Ur_1)$ have a uniform neighbourhood of zero not containing
non-trivial subgroups?
\end{q}

\begin{q} 
\bookproblem%
Is $\Iso(\Ur_1)$ with the uniform topology a Banach--Lie group?
\end{q}

The authors of \cite{RS} have established the following result
as an application of a new automatic continuity-type theorem and 
Example \ref{ex:homeo} above. 

\begin{theorem}[Rosendal and Solecki \cite{RS}]
The group $\Aut(\Q,\leq)$, considered as a discrete group, has the
{\em fixed point on metric compacta property,} that is, every action of
$\Aut(\Q,\leq)$ on a compact metric space by homeomorphisms has a common fixed point.
The same is true of the group $\Homeo_+[0,1]$.
\end{theorem}

This is particularly surprising in view of the Veech theorem, or, rather, its earlier version established by Ellis \cite{E60}: every discrete group $G$ acts freely on a suitable compact space by homeomorphisms (e.g. on $\beta G$). The two results seem to nearly contradict each other!

\begin{q}
\bookproblem%
Does the unitary group $U(\ell^2)$, viewed as a discrete group, have the 
fixed point on metric compacta property?
\end{q}

\begin{q}
\bookproblem%
The same question for the isometry group of the Urysohn space $\Ur_1$
of diameter one.
\end{q}

Extreme amenability is a strong form of {\em amenability,} an important
classical property of topological groups. A topological group $G$ is 
{\em amenable} if every compact
$G$-space admits an invariant probability Borel measure. Another reformulation:
the space $\RUCB(G)$ of all bounded right uniformly continuous real-valued
functions on $G$ admits a left-invariant mean, that is, a positive functional
$\phi$ of norm $1$ and the property $\phi(^gf) = \phi(f)$ for all
$g\in G$, $f\in\RUCB(G)$, where $^gf(x)=f(g^{-1}x)$. (Recall that 
the {\em right uniform structure} on $G$ is generated by
entourages of the diagonal of the form $V_R=\{(x,y)\in G\times G\colon
xy^{-1}\in V\}$, where $V$ is a neighbourhood of identity. For the
{\em left} uniformity, the formula becomes $x^{-1}y\in V$.)
For a general reference to amenability, see e.g. \cite{Pat}.

\begin{q}[A. Carey and H. Grundling \cite{CG}]
\bookproblem%
Let $X$ be a smooth compact manifold, and let $G$ be a
compact (simple) Lie group. Is the group $C^\infty(X,G)$ of all smooth
maps from $X$ to $G$, equipped with the pointwise operations and
the $C^\infty$ topology, amenable? 
\end{q}

This question is motivated by gauge theory models of 
mathematical physics \cite{CG}.

\begin{q}
\bookproblem%
To begin with, is the group of all continuous maps $C([0,1],SO(3))$ with
the topology of uniform convergence amenable?
\end{q}

The following way to prove extreme amenability of topological groups was
developed by Gromov and Milman \cite{GrM}. 
A topological group $G$ is called a
{\em L\'evy group} if there exists an increasing net $(K_\alpha)$ of compact
subgroups whose union is everywhere dense in $G$, having the following
property. Let $\mu_\alpha$ denote the Haar measure on the group $K_\alpha$,
normalized to one ($\mu_\alpha(K_\alpha)=1$). If $A\subseteq G$ is a Borel
subset such that $\liminf_\alpha\mu_\alpha(A\cap K_\alpha)>0$, then for every
neighbourhood $V$ of identity in $G$ one has
$\lim_\alpha \mu_\alpha(VA\cap K_\alpha)=1$. (Such a family of compact subgroups
is called a {\em L\'evy family.})

\begin{theorem}[Gromov and Milman \cite{GrM}] 
\label{th:gm}
Every L\'evy group is extremely amenable.
\end{theorem}

\begin{proof} We will give a proof in the case of a second-countable 
$G$, where one can assume the net $(K_\alpha)$ to be an increasing sequence. 
For every free ultrafilter $\xi$ on $\N$ the formula 
$\mu(A) =\lim_{n\to\xi}\mu_n(A\cap K_n)$ defines a {\em finitely-additive}
measure on $G$ of total mass one, invariant under multiplication on the left
by elements of the everywhere dense subgroup 
$\underline{G}=\cup_{n=1}^\infty K_n$. Besides, $\mu$ has the
property that if $\mu(A)>0$, then for every non-empty open $V$ one has
$\mu(VA)=1$. 
Let now $G$ act continuously on a 
compact space $X$. Choose an arbitrary $x_0\in X$. The push-forward, $\nu$, of the
measure $\mu$ to $X$ along the corresponding orbit map, given by
$\nu(B) = \mu\{g\in G\colon gx_0\in B\}$, is again a finitely-additive Borel
measure on $X$ of total mass one, invariant under translations by $\underline{G}$
and having the same ``blowing-up'' property: if $\nu(B)>0$ and $V$ is a
non-empty open subset of $G$, then $\nu(VB)=1$. Given a finite cover $\gamma$
of $X$, an element $W$ of the unique uniformity on $X$, and a finite subset
$F$ of $\underline{G}$, there is at least one $A\in\gamma$ with $\nu(A)>0$, consequently
$\nu(W[A])=1$ and for all $g\in F$ the translates $g\cdot W[A]$, having full
measure each, must overlap. This can be used to 
construct a Cauchy filter $\mathcal F$ of closed subsets of $X$ 
with $A\in{\mathcal F}$,
$g\in \underline{G}$ implying $gA\in{\mathcal F}$. The only point of 
$\cap{\mathcal F}$ is fixed under the action of $\underline{G}$ and therefore
of $G$ as well.
\end{proof}

For instance, the groups in Examples \ref{ex:one} and \ref{ex:gfw} are L\'evy
groups, and so is the isometry group $\Iso(\Ur)$ with the Polish topology
\cite{P05b}.

Theorem of Gromov and Milman cannot be inverted, because the extremely
amenable groups from Example \ref{ex:homeo} are not L\'evy: they simply do not
contain any non-trivial compact subgroups. What if such subgroups are present?
The following is a reasonable general reading of an old 
question by Furstenberg discussed at the end of \cite{GrM}.

\begin{q}
\bookproblem%
\label{q:furst}
Suppose $G$ is an extremely amenable
topological group containing a net of compact subgroups $(K_\alpha)$ whose
union is everywhere dense in $G$. Is $G$ a L\'evy group?
\footnote{I. Farah and S.Solecki have announced a counter-example (May 2006).}
\end{q}

\begin{q}
\bookproblem%
\label{q:furst2}
Provided the answer is yes, is the family $(K_\alpha)$ a
L\'evy family?\footnote{Cf. the previous footnote.}
\end{q}

A candidate for a ``natural'' counter-example is the group $SU(\infty)$, the
inductive limit of the family of special unitary groups of
finite rank embedded one into the other via
$SU(n)\ni V\mapsto 
\left(\begin{matrix}{\mathbb I} & 0 \\ 0 & V\end{matrix}\right)\in SU(n+1)$.
Equip $SU(\infty)$ with the inductive limit topology, that is, the finest
topology inducing the given topology on each $SU(n)$. 

\begin{q}
\bookproblem%
Is the group $SU(\infty)$ with the inductive limit topology extremely
amenable?
\end{q}

If the answer is yes, then Questions \ref{q:furst} and \ref{q:komega}
are both answered in the negative.

Historically the first example of an extremely amenable group was constructed
by Herer and Christensen \cite{HC}. Theirs was an abelian topological group
without strongly continuous unitary representations in Hilbert spaces
(an {\em exotic group}).

\begin{q}
\bookproblem%
Is the exotic group constructed in \cite{HC} a L\'evy group?
\end{q}

The following result shows that the properties of L\'evy groups are
diametrically opposed to those of locally compact groups in the setting of
ergodic theory as well as topological dynamics.

\begin{theorem}[Glasner--Tsirelson--Weiss \cite{GTW}]
Let a Polish L\'evy group act in a Borel measurable way on a Polish space $X$.
Let $\mu$ be a Borel probability measure on $X$ invariant under the action of
$G$. Then $\mu$ is supported on the set of $G$-fixed points.
\end{theorem}

\begin{q}[Glasner--Tsirelson--Weiss, {\em ibid.}] 
\bookproblem%
Is the same conclusion true if one 
only assumes
that the measure $\mu$ is {\em quasi-invariant}
under the action of $G$, that is,
for all $g\in G$ and every null-set $A\subseteq X$, the set $gA$ is null?
\end{q}

Recall that a compact $G$-space $X$ is called {\em minimal} if the orbit of
every point is everywhere dense in $X$. To every topological group $G$ 
there is associated the {\em universal minimal flow,} $\M(G)$, which is
a minimal compact $G$-space
uniquely determined by the property that every other 
minimal $G$-space is an image
of $\M(G)$ under an equivariant continuous surjection. (See \cite{Aus}.)
For example, $G$ is extremely amenable if and only if $\M(G)$ is a
singleton. If $G$ is compact, then $\M(G)=G$, but
for locally compact non-compact groups, starting with $\Z$,
the flow $\M(G)$ is typically very
complicated and highly non-constructive, in particular it is never metrizable
\cite{KPT}. A discovery of the recent years has been that even non-trivial
universal minimal flows of 
``large'' topological groups are sometimes manageable.

\begin{example} The flow $\M(\Homeo_+(\s^1))$ is the circle $\s^1$ itself,
equipped with the canonical action of the group $\Homeo_+(\s^1)$ of
orientation-preserving homeomorphisms, with the compact-open topology
\cite{P98a}.
\end{example}

\begin{example} Let $S_\infty$ denote the infinite symmetric group, that is,
the Polish group of all bijections of the countably 
infinite discrete space $\omega$ onto itself, equipped with the topology of
pointwise convergence. The flow $\M(S_\infty)$ can be identified with the
set of all linear orders on $\omega$ with the topology induced from
$\{0,1\}^{\omega\times\omega}$ under the identification of each order with
the characteristic function of the corresponding relation \cite{Gl-W1}.
\end{example}

\begin{example} Let $C=\{0,1\}^\omega$ stand for the Cantor set. The
minimal flow $\M(\Homeo(C))$ can be identified with the space of
all maximal chains of closed subsets of $C$, equipped with the Vietoris
topology. This is the result of Glasner and Weiss \cite{Gl-W2}, while the
space of maximal chains was introduced into the dynamical context by
Uspenskij \cite{Usp00}.
\end{example}

\begin{q}[Uspenskij]
\bookproblem%
Give an explicit description of the universal minimal flow
of the homeomorphism group $\Homeo(X)$ of a closed compact manifold $X$
in dimension $\dim X>1$ (with the compact-open topology). 
\label{q:umx}
\end{q}

\begin{q}[Uspenskij]
\bookproblem%
The same question for the group of homeomorphisms of the Hilbert
cube $Q=\I^\omega$.
\label{q:umq}
\end{q}

Note that both $X$ and $Q$ form minimal flows for the respective
homeomorphism groups, but they are not universal \cite{Usp00}.
Interesting recent advances on both questions \ref{q:umx} and \ref{q:umq} 
belong to Yonatan Gutman \cite{Gut}.

\begin{q}[Uspenskij]
\bookproblem%
Is the pseudoarc $P$ the universal minimal flow for its own homeomorphism group?
\end{q}

A recent investigation \cite{IS} might provide means to
attack this problem.

Let $G$ be a topological group. 
The completion of $G$ with regard to the left
uniform structure (the {\em left} completion), denoted by $\hat G^L$, is 
a topological semigroup with jointly continuous multiplication
\cite[Prop. 10.2(a)]{RD},
but in general not a topological group \cite{Di}. 
Note that every left uniformly continuous
real-valued function $f$ on $G$ extends to a unique continuous function
$\hat f$ on $\hat G^L$. 
Say that such an $f$ is {\em oscillation stable} if for every 
$\e>0$ there is a right ideal $J$ in the topological semigroup $\hat G^L$
with the property that the values of $\hat f$ at any two points of $J$ differ 
by $< \e$. 
If $H$ is a closed subgroup of $G$, say that the homogeneous space $G/H$
is {\em oscillation stable} if every bounded left uniformly continuous 
function $f$ on $G$ that factors through the quotient map
$G\to G/H$ is oscillation stable. If $G/H$ is not oscillation stable, we say
that $G/H$ has {\em distortion}.

\begin{example}
The unit sphere $\s^\infty$ in the separable Hilbert space $\ell^2$,
considered as the homogeneous factor-space of the unitary group $U(\ell^2)$
with the strong topology, has distortion. It means that there exists a
uniformly continuous function $f\colon\s^\infty\to\R$ whose range of values
on the intersection of $\s^\infty$ with every infinite-dimensional linear
subspace contains the interval (say) $[0,1]$.
This is a famous and very difficult
result by Odell and Schlumprecht \cite{OS}, answering a 30 year-old problem.
The following question is well-known in geometric functional analysis.
\end{example} 

\begin{q}
\bookproblem%
Does there exist a direct proof of Odell and Schlumprecht's result, based on
the intrinsic geometry of the unit sphere and/or the unitary group?
\label{q:dirdist}
\end{q}

\begin{example} 
The set $[\Q]^n$ of all $n$-subsets of $\Q$,
considered as a homogeneous factor-space of $\Aut(\Q,\leq)$, is oscillation
stable if and only if $n=1$. For $n=1$, oscillation stability simply means that
for every finite colouring of $\Q$, there is a monochromatic subset $A$
order-isomorphic to $\Q$ (this is obvious). For $n\geq 2$, distortion
of $[\Q]^n$ means the existence of a finite colouring of this set with
$k\geq 2$ colours
such that for every subset $A$ order-isomorphic to the rationals the set
$[A]^n$ contains points of all $k$ colours. This follows easily from classical
Sierpi\'nski's partition argument \cite{Sie1}, cf. 
\cite[Example 5.1.27]{P05}.
\end{example}

The above setting for analysing distortion/oscillation stability in the
context of topological transformation groups was proposed in \cite{KPT} and
discussed in \cite{P05}.
The most substantial general result within this framework is presently the following.

\begin{theorem}[Hjorth \cite{hjorth}]
Let $G$ be a Polish topological group. Considered as a $G$-space
with regard to the action on itself by left translations, $G$ has
distortion whenever $G\neq\{e\}$.
\end{theorem}

\begin{q}[Hjorth \cite{hjorth}]
\bookproblem%
Let $E$ be a separable Banach space and let $\s_E$ denote
the unit sphere of $E$ viewed as an $\Iso(E)$-space, where the latter group
is equipped with the strong operator topology. Is it true that 
the $\Iso(E)$-space $\s_E$ has distortion?
\end{q}

Note of caution: this would not, in general, mean that $E$ 
has distortion in the sense of theory of
Banach spaces \cite[Chapter 13]{BL}, as the two concepts only coincide
for Hilbert spaces.

For an ultrahomogeneous separable metric space $X$, oscillation stability of $X$ 
equipped with the standard action of the Polish group of isometries $\Iso(X)$
is equivalent to the following property. For every finite cover $\gamma$ of
$X$, there is an $A\in\gamma$ such that for each $\e>0$, the 
$\e$-neighbourhood of $A$ contains an isometric copy of $X$. 
The following could provide a helpful insight into question
\ref{q:dirdist}.

\begin{q}
\bookproblem%
Is the metric space $\Ur_1$ oscillation stable? 
\label{q:uos}
\end{q}

The Urysohn metric space $\Ur$ itself has distortion, but for trivial reasons,
just like any other
unbounded connected ultrahomogeneous metric space.

The oscillation stability of a metric space $X$ whose distance assumes 
a discrete collection of values is equivalent to the property that whenever
$X$ is partitioned into two subsets, at least one of them contains an
isometric copy of $X$.
The Urysohn metric space $\Ur_{\{0,1,2\}}$ universal for the class of metric
spaces whose distances take values $0,1,2$ is oscillation
stable, because it is isometric to the path metric space associated to the
{\em infinite random graph} $R$, and oscillation stability is an immediate
consequence of an easily proved property of $R$ known as {\em indestructibility}
(cf. \cite{Cam}).
Very recently, Delhomme, Laflamme, Pouzet, and Sauer \cite{DLPS} have established
oscillation stability of the universal Urysohn metric space
$\Ur_{\{0,1,2,3\}}$ with the distance taking values $0,1,2,3$.
The following remains unknown.

\begin{q}
\bookproblem%
Let $n\in\N$, $n\geq 4$. Is the  
universal Urysohn metric space $\Ur_{\{0,1,\ldots,n\}}$ oscillation stable?
\end{q}

Resolving the following old question may help.

\begin{q}[M. Fr\'echet \cite{Fre}, p. 100; P.S. Alexandroff \cite{U72}]
\bookproblem%
Find a model for the Urysohn space $\Ur$, that is, a
concrete realization.
\end{q}

Several such models are known for the random graph (thence,
$\Ur_{\{0,1,2\}}$), cf. \cite{Cam}.
\begin{q}
\bookproblem%
Find a model for the metric space $\Ur_{\{0,1,2,3\}}$.
\end{q}

In connection with Uspenskij's examples of universal 
second-countable topological groups \cite{Usp86,Usp90}, including 
$\Iso(\Ur)$, the following remains unresolved.

\begin{q}[V.V. Uspenskij \cite{Usp98}]
\bookproblem%
Does there exist a universal topological group of every given
infinite weight $\tau$? 
\end{q}

\begin{q}[V.V. Uspenskij]
\bookproblem%
The same, for {\it any} uncountable weight?
\end{q}

\begin{q}[A.S. Kechris]
\bookproblem%
Does there exist a {\em co-universal} Polish
topological group $G$, that is, such that every other Polish group is
a topological factor-group of $G$?
\end{q}

In the abelian case, the answer is in the positive \cite{SPW}.

\begin{q}[A.S. Kechris]
\bookproblem%
Is every Polish topological group a 
topological factor-group of a subgroup of $U(\ell^2)$ with the strong
topology?
\label{q:fug}
\end{q}

Again, in the abelian case the answer is in the positive \cite{GaP}.

\begin{q}
\bookproblem%
Is the free topological group $F(X)$  \cite{graev} on a metrizable compact space
isomorphic to a topological subgroup of the unitary group $U({\mathcal H})$
of a suitable Hilbert space $\mathcal H$, equipped with the strong topology?
\label{q:free}
\end{q}

Galindo has announced \cite{galindo} a
positive answer for free abelian topological groups. Uspenskij \cite{Usp06} 
has given a very elegant proof of a more general result: the free abelian topological group $A(X)$ of a Tychonoff space embeds into $U({\mathcal H})$ as a topological subgroup. This suggests a more general vesion of the same question:

\begin{q}
\bookproblem%
The same question for an arbitrary Tychonoff space $X$.
\label{q:zeroone}
\end{q}

In connection with questions \ref{q:fug}, \ref{q:free} and \ref{q:zeroone},
let us remind the following old problem.

\begin{q}[A.I. Shtern \cite{Shtern0}]
\bookproblem%
What is the intrinsic characterization of topological subgroups of 
$U(\ell^2)$ (with the strong topology)?
\label{q:schtern}
\end{q}

A unitary representation $\pi$ of a topological group $G$ in a Hilbert space
$\H$ (that is, a strongly continuous homomorphism $G\to U(\H)$)
{\it almost has invariant vectors}
if for every compact $F\subseteq G$ and every $\e>0$ there is a
$\xi\in\H$ with $\norm\xi=1$
and $\norm{\pi_g\xi-\xi}<\e$ for every $g\in F$. 
A topological group $G$ has
{\em Kazhdan's property (T)} if, whenever a
unitary representation of $G$
almost has invariant vectors, it has an invariant vector of norm one.
For an excellent account of this rich theory, see the book
\cite{dlHV} and especially its many times
extended and updated English version, currently in preparation and
available on-line \cite{BdlHV}. 

Most of the theory is concentrated in the locally compact case. Bekka has
shown in \cite{Bek3} that the group $U(\ell^2)$ with the strong topology
has property $(T)$.

\begin{q}[Bekka \cite{Bek3}]
\bookproblem%
Does the group $U(\ell^2)$ with the uniform topology have property $(T)$?
\end{q}

\begin{q}[Bekka \cite{Bek3}]
\bookproblem%
Does the unitary group $U(\ell^2(\Gamma))$ of a non-separable Hilbert space
($\abs\Gamma>\aleph_0$),
equipped with the strong topology, have property $(T)$?
\end{q}

Here is a remarkable ``large'' topological group that has been receiving much attention recently. Let $\norm\cdot_2$ denote the 
Hilbert-Schmidt norm on the $n\times n$
matrices, $\norm{A}_2=\left(\sum_{i,j=1}^n\abs{a_{ij}}^2\right)^{1/2}$, and
let $d_n$ be the {\em normalized} Hilbert-Schmidt metric
on the unitary group $U(n)$, that is, $d_n(u,v) = \frac 1{\sqrt n}\norm{u-v}_2$.
Choose a free ultrafilter $\xi$ on the natural numbers and
denote by $U(\xi)_2$ the factor-group of the direct product $\prod_{n\in\N} U(n)$
by the normal subgroup 
$N_\xi = \{(x_n)\colon \lim_{n\to\xi}d_n(e,x_n)=0\}$.

The following question is a particular case of Connes' Embedding Conjecture
\cite{Con}, for a thorough discussion 
see \cite{O} and references therein.

\begin{q}[Connes' Embedding Conjecture for Groups]
\bookproblem%
Is every countable group isomorphic to a subgroup of $U(\xi)_2$ (as an
abstract group)? 
\end{q}

Groups isomorphic to subgroups of $U(\xi)_2$ are called {\em hyperlinear.} 
Here are some of the most important particular cases of the above problem.

\begin{q}
\bookproblem[3]
Are countable groups from the following classes hyperlinear:
(a) one-relator groups; (b) hyperbolic groups \cite{Gr0}; (c) groups amenable at infinity (a.k.a. topologically amenable groups, exact groups) \cite{ADR}?
\end{q}

Under the natural
bi-invariant metric $d(x,y) =\lim_{n\to\xi} d_n(x_n,y_n)$, 
the group $U(\xi)_2$ is a complete non-separable metric group whose left and
right uniformities coincide, isomorphic to a topological subgroup
of $U(\ell^2({\mathfrak c}))$ with the strong topology. 
Understanding the topological group structure of $U(\xi)_2$ may prove important.

The Connes' Embedding Conjecture itself can be reformulated in the language of topological groups as follows. Say, following \cite{PU}, that a topological group $G$ has {\em Kirchberg's property} if, whenever $A$ and $B$ are finite subsets of $G$ with the property that every elemant of $A$ commutes
with every element of $B$, there exist finite subsets $A'$ and $B'$ of $G$
that are arbitrarily close to $A$ and $B$, respectively,
such that every element of $A'$
commutes with every element of $B'$, and the subgroups of $G$ generated by $A'$ and $B'$ are relatively compact. As noted in \cite{PU}, the deep results of \cite{kirchberg}, modulo a criterion from \cite{el}, immediately imply that the Connes Embedding Conjecture is equivalent to the statement that the unitary group $U(\ell^2)$ with the strong topology has Kirchberg's property.

\begin{q}
\bookproblem[2]%
Do the following topological groups have Kirchberg's property: (a) the infinite symmetric group $S_\infty$, (b) the group $\Aut(X,\mu)$ of measure-preserving transformations of a standard Lebesgue measure space with the coarse topology?
\end{q}

It was shown in \cite{PU} that $\Iso(\Ur)$ has Kircherg's property.


\begin{thebibliography}{10}
\bibitem{AHK}
M. Ajtai, I. Havas and J. Koml\'os, 
\newblock{Every group admits a bad topology.}
\newblock
In: {\em Studies in pure. mathematics}, 21-34 (P. Erdos, ed.) 
Birkhauser, Basel, 1983.

\bibitem{ADR}  C. Anantharaman-Delaroche, J. Renault.
\newblock
Amenable groupoids. 
\newblock
Monographies de L'Enseignement Math\'ematique, 36,
L'Enseignement MathŽmatique, Geneva, 2000.

\bibitem{Aus} J. Auslander,
\newblock{Minimal Flows and Their Extensions.}
\newblock North-Holland Mathematics Studies 153,
North-Holland, Amsterdam--NY--London--Tokyo, 1988.

\bibitem{Ba} W. Banaszczyk,
\newblock {On the existence of exotic Banach--Lie groups.} 
\newblock {\em Math. Ann.} 264:485--49, 1983.

\bibitem{Bek3} M.B. Bekka, 
\newblock
Kazhdan's property (T) for the unitary group of a separable 
Hilbert space. \newblock
{\em Geom. Funct. Anal.} 13:509--520, 2003.

\bibitem{BdlHV} B. Bekka, P. de la Harpe, and A. Valette,
\newblock
{Kazhdan's Property $(T)$.} \newblock 
Book in preparation, current version available
from {\tt http://poncelet.sciences.univ-metz.fr/$\sim$bekka/}

\bibitem{BL} Y. Benyamini and J. Lindenstrauss,
\newblock{Geometric Nonlinear Functional Analysis,} Vol. 1.
\newblock
Colloquium Publications 48, American Mathematical Society,
Providence, RI, 2000.

\bibitem{Cam} P. Cameron, 
\newblock{The random graph.}
\newblock In:
The Mathematics of Paul Erdos, 331-351, J. Ne\u setril, R. L. Graham, eds.,
Springer, 1996.

\bibitem{CG} A. Carey and H. Grundling, 
\newblock{On the problem of the amenability of the
gauge group.}\newblock
{\em Lett. Math. Phys.} 68:113--120, 2004.

\bibitem{Con} A. Connes, \newblock
{Classification of injective factors}. \newblock
{\em Ann. of Math.} 104:73--115, 1976.

\bibitem{CotRic} M. Cotlar and R. Ricabarra,
\newblock On the existence of characters in topological groups.
\newblock {\em Amer. J. Math.} 76:375--388, 1954.

\bibitem{DLPS} C. Delhomme, C. Laflamme, M. Pouzet, and N. Sauer,
\newblock Divisibility of countable metric spaces.
\newblock ArXiv e-print math.CO/0510254.

\bibitem{Di} J. Dieudonn\'e, 
\newblock{Sur la completion des groupes topologiques.}
\newblock
C. R. Acad. Sci. Paris 218:774--776, 1944.

\bibitem{DPS} D. Dikranjan, I. Prodanov and L. Stoyanov,  
\newblock{Topological groups: characters, dualities and minimal group topologies.}
\newblock {\em
Monographs and Textbooks in Pure and Applied Mathematics} 130, 
Marcel Dekker Inc., New York-Basel, 1990.

\bibitem{E60} 
R. Ellis, Universal minimal sets, {\it Proc. Amer. Math. Soc.} {\bf
11} (1960), 540--543.

\bibitem{EK} R. Ellis and H.B. Keynes,
\newblock
{Bohr compactifications and a result of F\o lner.}
\newblock {\em Israel J. Math.} 12:314--330, 1972.

\bibitem{el}  R. Exel and T.A. Loring, 
\newblock{Finite-dimensional representations of free product 
$C\sp *$-algebras.} \newblock Internat. J. Math. 3:469--476, 1992.

\bibitem{Fol} E. F\o lner, 
\newblock Generalization of a
theorem of Bogoliuboff to topological abelian groups.
\newblock {\em Math. Scand.} 2:5--18, 1954.

\bibitem{Fre} M. Fr\'echet, \newblock{Les espaces abstraits.} 
\newblock Paris, 1928.

\bibitem{galindo} J. Galindo, \newblock
On unitary representability of topological groups. \newblock
Preprint, 2005.

\bibitem{GK} S. Gao and A.S. Kechris, 
\newblock{On the Classification of 
Polish Metric Spaces up to Isometry.} \newblock 
Memoirs of the Amer. Math. Soc. 766, 2003.

\bibitem{GaP} S. Gao and V. Pestov,
\newblock{On a universality property of some abelian Polish groups.}
\newblock {\em
Fund. Math.}  179:1--15, 2003.

\bibitem{GP2} 
T. Giordano and V. Pestov, 
\newblock Some extremely amenable groups related to
operator algebras and ergodic theory.
\newblock ArXiv e-print math.OA/0405288, to appear in J. Inst. Math. Jussieu.

\bibitem{Gl1} S. Glasner,
\newblock
On minimal actions of Polish groups,
\newblock {\em Top. Appl.} 85:119--125, 1998.

\bibitem{GTW} E. Glasner, B. Tsirelson, and B. Weiss,
\newblock{The automorphism group of the Gaussian measure cannot act pointwise.}
\newblock ArXiv e-print math.DS/0311450. To appear in
{\em Israel J. Math.}

\bibitem{Gl-W1}
E. Glasner and B. Weiss, 
\newblock{Minimal actions of the group $S(\mathbb Z )$ of
permutations of the integers}. \newblock
{\em Geom. and Funct. Anal.} 12:964--988, 2002.

\bibitem{Gl-W2}
E. Glasner and B. Weiss, 
\newblock{The universal minimal system for the group
of homeomorphisms of the Cantor set.}
\newblock {\em Fund. Math.} 176:277--289, 2003.

\bibitem{graev}
M.I. Graev, \newblock
Free topological groups. \newblock
Amer. Math. Soc. Translation 1951, 35,  61 pp., 1951.

\bibitem{Gr0} M. Gromov,
\newblock Hyperbolic groups. 
\newblock {\em
Essays in group theory,} 75--263, 
Math. Sci. Res. Inst. Publ., 8, 
Springer, New York, 1987. 

\bibitem{Gr} M. Gromov, 
\newblock{Metric Structures for Riemannian and
Non-Riemannian Spaces.} 
\newblock {\em Progress in Mathematics} 152, Birkhauser
Verlag, 1999.

\bibitem{GrM}
M. Gromov and V.D. Milman,
\newblock A topological application of the isoperimetric inequality.
\newblock {\em
Amer. J. Math.}, 105: 843--854, 1983.

\bibitem{Gut} Y. Gutman, \newblock
Minimal actions of homeomorphism groups, \newblock
preprint, 2005.

\bibitem{dlHV} P. de la Harpe and A. Valette, \newblock
{La propri\'et\'e {\rm (T)} de Kazhdan pour les groupes
localements compacts,}
\newblock {\em Ast\'erisque} 175, 1989.

\bibitem{HC}  W. Herer and J.P.R. Christensen,
\newblock On the existence of pathological
submeasures and the construction of exotic topological groups.
\newblock {\em Math. Ann.} 213:203--210, 1975.

\bibitem{hjorth} G. Hjorth, 
\newblock{An oscillation theorem for groups of
isometries.}
\newblock preprint, Dec. 31, 2004, 28 pp.

\bibitem{IS} T. Irwin and S. Solecki, \newblock
Projective Fra\"\i ss\'e limits and the pseudo-arc.
\newblock Trans. Amer. Math. Soc., to appear.

\bibitem{KPT} A.S. Kechris, V.G. Pestov and S. Todorcevic,
\newblock Fra\"\i ss\'e limits, Ramsey theory, and topological dynamics of
automorphism groups.
\newblock {\em Geom. Funct. Anal.} 15:106--189, 2005.

\bibitem{kirchberg}  E. Kirchberg, \newblock{On nonsemisplit extensions, 
tensor products and exactness of group $C\sp *$-algebras.}
\newblock Invent. Math. 112:449--489, 1993.

\bibitem{OS}
E. Odell and T. Schlumprecht, 
\newblock{The distortion problem.}
\newblock {\em Acta Math.} 173:259--281, 1994.

\bibitem{O}  N. Ozawa, 
\newblock{About the QWEP conjecture.}
\newblock {\em
Internat. J. Math.} 15:501--530, 2004. 

\bibitem{Pat} A.T. Paterson,
\newblock
{Amenability.} \newblock Math. Surveys and Monographs
29, Amer. Math. Soc., Providence, RI, 1988.

\bibitem{P98a} V.G. Pestov, 
\newblock
On free actions, minimal flows, and a problem by Ellis,
\newblock {\em
Trans. Amer. Math. Soc.} 350:4149-4165, 1998.

\bibitem{P98b} 
V. Pestov, \newblock{Some universal constructions in abstract 
topological dynamics.} \newblock
{\em Contemporary Math.} 215:83--99, 1998.

\bibitem{P02} V. Pestov, 
\newblock{Ramsey--Milman phenomenon, 
Urysohn metric spaces,
and extremely amenable groups.} \newblock
{\em Israel Journal of Mathematics}
127:317-358, 2002. \newblock{Corrigendum,} {\em ibid.} 145:375-379, 2005.

\bibitem{P05} V. Pestov, 
\newblock{Dynamics of infinite-dimensional groups and 
Ramsey-type phenomena.} 
\newblock {\em Publica\c c\~oes dos Col\'oquios de Matem\'atica,}
IMPA, Rio de Janeiro, 2005.

\bibitem{P05b} V. Pestov, \newblock
The isometry group of the Urysohn space as a L\'evy group.
\newblock  To appear in: Proceedings of the 6-th Iberoamerican Conference on Topology and its Applications (Puebla, Mexico, 4-7 July 2005).
\newblock ArXiv e-print math.GN/0509402.

\bibitem{PU} V.G. Pestov and V.V. Uspenskij,
\newblock{Representations of residually finite groups by isometries of the Urysohn space.} \newblock
ArXiv e-print math.RT/0601700, to appear in J. Ramanujan Math. Soc.

\bibitem{PZ} I. Protasov, Y. Zelenyuk, 
\newblock{Topologies on Groups Determined 
by Sequences.} 
\newblock Lviv, VNTL Publishers, 1999.

\bibitem{RD} W. Roelcke and S. Dierolf, 
\newblock{Uniform Structures on Topological
Groups and Their Quotients.} \newblock McGraw-Hill, 1981. 

\bibitem{RS} C. Rosendal and S. Solecki,
\newblock{Automatic continuity of group
homomorphisms and discrete groups with the fixed point on metric
compacta property}. \newblock  preprint, 2005.

\bibitem{Shtern0} A. Shtern,
{\it Unitary Representation,} in: {\it Mathematical Encyclopaedia,}
Vol. 5, Sov. Encycl., 1984, pp. 508--513 (in Russian).

\bibitem{SPW}
D. Shakhmatov, J. Pelant, and S. Watson, 
\newblock{A universal complete metric abelian group of a given weight.}
\newblock {\em
Bolyai Soc. Math. Stud.} 4:431--439, 1995.

\bibitem{Sie1} W. Sierpi\'nski, 
\newblock{Sur un probl\`eme de la th\'eorie des relations.}
\newblock {\em Ann. Scuola Norm. Sup. Pisa} 2:285--287, 1933.

\bibitem{U25}
P. S. Urysohn, 
\newblock{Sur un espace metrique universel.}
\newblock {\em
C. R. Acad. Sci. Paris} 180:803--806, 1925.

\bibitem{U72} P.S. Urysohn, \newblock{On the universal metric space.}
\newblock
In: P.S. Urysohn. Selected Works, vol. 2, 747--769, P. S. Alexandrov, ed., 
Nauka, Moscow, 1972 (in Russian).

\bibitem{Usp86} V.V. Uspenskij,
\newblock
{A universal topological group with countable base.}
\newblock {\em
Funct. Anal. Appl.} 20:160--161, 1986.

\bibitem{Usp90}
V.V. Uspenski\u\i,  
\newblock On the group of isometries of the
Urysohn universal metric space.
\newblock {\em Comment. Math. Univ. Carolinae} 31:181-182, 1990.

\bibitem{Usp98} V.V. Uspenskij, 
\newblock On subgroups of minimal
topological groups.
\newblock ArXiv e-print math.GN/0004119.

\bibitem{Usp00} 
V. Uspenskij, 
\newblock{On universal minimal compact $G$-spaces.}
\newblock Topology Proc. 25:301--308, 2000.

\bibitem{Usp06} 
V.V. Uspenskij, 
\newblock{Unitary representability of free abelian topological groups.}
\newblock ArXiv e-print math.RT/0604253.

\bibitem{V} W.A. Veech, \newblock{Topological dynamics.}
\newblock {\em Bull. Amer. Math. Soc.} 83:775--830, 1977.

\bibitem{V2} W.A. Veech, 
\newblock{The equicontinuous structure relation for 
minimal abelian transformation groups.} 
\newblock Amer. J. Math. 90:723--732, 1968.

\bibitem{Ver04} A.M. Vershik, 
\newblock{The universal and random metric
spaces.} 
\newblock {\em Russian Math. Surveys} 356:65--104, 2004.

\bibitem{W} B. Weiss, 
\newblock{Single Orbit Dynamics.}
\newblock 
CBMS Regional Conference Series in Mathematics 95, 
American Mathematical Society, Providence, RI, 2000.

\end{thebibliography}
\end{document}